\documentclass[10pt]{article}
\newtheorem{theorem}{Theorem}

\newtheorem{remark}{Remark}

\def\Frac#1#2{\frac{\displaystyle{#1}}{\displaystyle{#2}}}
\usepackage{amssymb}
\usepackage{epsfig}

\begin{document}
 \title{{\Large\bf On the complex zeros of Airy and Bessel functions and those of their derivatives}}

\author{
\\
A. Gil\\
Departamento de Matem\'atica Aplicada y CC. de la Computaci\'on.\\
Universidad de Cantabria. 39005-Santander, Spain.\\
\and
\\
J. Segura\footnote{Corresponding author}\\
        Departamento de Matem\'aticas, Estad\'{\i}stica y 
        Computaci\'on.\\
        Universidad de Cantabria, 39005-Santander, Spain.\\  
}

\date{\ }

\maketitle

\begin{center}
{\it\large Dedicated to the memory of Frank W. J. Olver}
\vspace*{1cm}
\end{center}

%\begin{abstract}
\begin{center}
{\bf Abstract}
\end{center}

\noindent
We study the distribution of zeros of general solutions of the Airy and Bessel 
equations in the complex plane. Our results characterize the 
patterns followed by the zeros for any solution, in such a way that if one zero
is known it is possible to determine the location of the rest of zeros. 
%\end{abstract}

\vspace*{0.2cm}


\section*{{\large Introduction}}

Airy and Bessel functions are important examples of special functions satisfying second order lineal
ODEs. Their zeros appear in a great number of applications in physics and engineering. Particularly, 
the complex zeros are quantities appearing in some problems of quantum physics \cite{Cruz:1982:ZHF}, electromagnetism \cite{Davies:1973:CZO} and wave propagation and scattering \cite{Parnes:1972:CZO,Ferreira:2008:ZOT}.

 Information is available regarding the distribution of zeros of some particular solutions of these
differential equations. For example, the zeros of Airy functions ${\rm Ai}(z)$, ${\rm  Bi}(z)$ and
$i m {\rm Ai}(z)-Bi(z)$ ($m$ real) are studied in \cite{Olver:1954:TAE} (see also \cite{Fab:1999:OTR})
together with the zeros of the Bessel functions $J_{n}(z)$ and $Y_{n}(z)$ of
integer order (see also \cite{Olver:2010:BES}) and the zeros of
semi-integer order for combinations $\cos\alpha J_{\nu}(z)-\sin\alpha Y_{\nu}(z)$ are discussed in 
\cite{Davies:1973:CZO}. The zeros of Hankel functions $H_{\nu}^{(1)}(z)$ and $H_{\nu}^{(2)}(z)$ for
real order $\nu$, which are also solutions of the Bessel equation, are analyzed in \cite{Cruz:1982:ZHF}.
However, up to date there is no available description of the distribution of zeros for general solutions
of the Airy and Bessel equations.

We perform this analysis by considering two complementary approaches. In the first place,
a qualitative picture of the possible patterns of the complex zeros is obtained by assuming that the Liouville-Green
approximation holds. This is the starting point of the numerical method of \cite{Segura:2012:CCZ}. In 
the second place, the use of asymptotics (both of Poincar\'e type and uniform asymptotics for large orders in the case
of Bessel functions) will provide more detailed and quantitative information. These results will characterize the possible
patterns followed by the zeros for any solution, leading to the development
of methods which are able to compute safely and accurately all the zeros in
a given region and having as only input data the location of just one zero (and
also the order for the case of the Bessel equation).

\section{{\large Liouville-Green approximation. Stokes and anti-Stokes lines.}}

An apparently naive method for studying the distribution of the complex zeros
of second order ODEs 
\begin{equation}
\label{ODE}
y''(z)+A(z)y(z)=0
\end{equation}
consist in considering that $A(z)$ is ``locally constant" in the sense we next describe.
 
We start by considering the trivial case of $A(z)$ constant. Then the
general solution of (\ref{ODE}) reads
$$
y(z)=C\sin \left(\sqrt{A(z)}\,(z-\psi)\right),
$$
and the zeros are over the line
$$
z=\psi+e^{-i\frac{\varphi}{2}}\lambda,\ \lambda\in{\mathbb R},\  \varphi=\arg{A(z)}.
$$
In other words, writing $z=u+iv$ we have that the zeros are over an integral line of
\begin{equation}
\label{LGA}
\Frac{dv}{du}=-\tan (\varphi/2 ).
\end{equation}
If $A(z)$ is constant the zeros of the derivative $y'(z)=0$ would lie on the same line (\ref{LGA}).
Of course, in general $A(z)$ will not be a constant and we consider the following simplifying
assumption: the curves where the zeros lie 
are also given by (\ref{LGA}), but with variable $\varphi$; this is what we mean when we say
that we consider that $A(z)$ is ``locally constant". For the zeros of $y'(z)$ the same 
approximation makes sense if the variation of $A(z)$ is sufficiently slow.

This assumption regarding the distribution of zeros is equivalent to considering
that the Liouville-Green approximation (LG) is accurate. The LG approximation (also called WKB or WKBJ approximation) 
with a zero at $z^{(0)}$ is
$$
y(z)\approx C A(z)^{-1/4}\sin\left(\int_{z^{(0)}}^z A(\zeta)^{1/2}d\zeta \right).
$$
Then, if $z^{(0)}$ is a zero, other zeros of the LG approximation lie over the curve such that
\begin{equation}
\label{ASL}
\Im \int_{z^{(0)}}^z A(\zeta)^{1/2}d\zeta=0,
\end{equation}
and those curves are also given by (\ref{LGA}). These are the so-called anti-Stokes lines
(ASLs) in the LG approximation. Similarly, the curve defined by 
\begin{equation}
\label{ASL}
\Re \int_{z^{(0)}}^z A(\zeta)^{1/2}d\zeta=0,
\end{equation}
gives the Stokes line (SL) passing through $z_0$.

It is a well-known fact that when SLs and ASLs intersect they do so perpendicularly. However, two ASLs (or SLs) 
can not intersect except at the zeros 
and poles of $A(z)$. 
This fact singularizes the role of the ASLs (or SLs) emerging from the zeros (also called turning points) and the poles. 
The ASLs (or SLs) are called principal when they emerge from the zeros of A(z). In particular, if $A(z)$ has a zero
of multiplicity $m$ at $z=z_0$ such that $A(z)=a(z-z_0)^m (1+{\cal O}(z-z_0))$ as $z\rightarrow z_0$, 
$m\in {\mathbb N}$, then $m+2$ principal ASLs emerge from $z_0$ in the directions
\begin{equation}
\arg (z-z_0) =\Frac{1}{m+2}(-\arg (a)+2k\pi),\,k=0,\ldots,m+1.
\end{equation}

Principal ASLs (SLs) divide the complex plane in different disjoint domains such that any ASL (SL) is either inside one of these domains or is a principal ASL (SL) itself. In the literature these domains are called anti-Stokes and Stokes domains.

As discussed in \cite{Segura:2012:CCZ}, the qualitative picture provided by the
ASLs gives generally a quite accurate picture of how the complex zeros may be distributed
in the complex plane.
For the zeros of $y'(z)$, as we will see, the approximation that they are on ASL curves given by (\ref{LGA}) 
(and equivalently \ref{ASL}) will be also in general sufficiently accurate. Furthermore, the numerical method
for complex zeros of $y(z)$ given in \cite{Segura:2012:CCZ} 
is also successful for the zeros of  $y'(z)$ with only a simple change in the 
iteration function of the fixed point method. 

\section{{\large Zeros of Airy functions}}

Next we analyze the distribution of the zeros of Airy functions.  We begin with 
a description of the ASLs and SLs, followed by an analysis of the distribution of
the zeros for the general solution $\cos\alpha {\rm Ai}(z)+\sin\alpha {\rm Bi}(z)$; this
analysis will be important in the description of the zeros of general Bessel functions
$\cos\alpha J_{\nu}(z)-\sin\alpha Y_{\nu}(z)$.

\subsection{{\normalsize Anti-Stokes lines for the Airy equation in the LG approximation}}
\label{ASLA}

For the Airy equation $y''(z)-A(z)y(z)=0$, with $A(z)=-z$, the ASL passing through a
point $z_i$ is given in the LG approximation by
\begin{equation}
\Im \int_{z_i}^{z}\sqrt{-\eta}d\eta =0.
\end{equation}
This is the curve $\Re (z^{3/2}-z_i^{3/2})=0$ ($z=r e^{i\theta}$, 
$z_i=r_i e^{i\theta_i}$). 

For $z_i=0$ we obtain the three principal ASLs emerging from the origin in the directions $\arg z = -\pi, \pm \pi/3$. 
The remaining ASLs (not having $z_i=0$ in the path) have the expression
\begin{equation}
\label{LGAir}
r(\theta) = r_0 \left|\cos \Frac{3\theta}{2} \right|^{-2/3}
\end{equation}
with $\theta$ running in three possible different intervals, and therefore with ASLs given by $r(\theta)e^{i\theta}$
in three different sectors 
$$S_j=\left\{z : \arg z \in \left((2 j-1)\Frac{\pi}{3},(2 j+1)\Frac{\pi}{3}\right)\right\},\, j=-1,0,1.$$ 
In each sector, the curve (\ref{LGAir})  approaches
asymptotically the boundaries of the sector as $\theta$ approaches the boundaries of the intervals, and it passes 
through the point $z=r_0 e^{2\pi i j/3}$.

The curves (\ref{LGAir}) can be related to lines parallel to the positive real axis by a simple transformations. 
This relation will be interesting when we later obtain asymptotic approximations. We define
\begin{equation}
z (\xi)=\left(\frac32 \xi \right)^{2/3},
\end{equation}
where the factor $3/2$ inside the parenthesis is introduced for later convenience. Then it is straightforward to verify
that the values $e^{i\pi/3}z (\xi)$ with
\begin{equation}
\xi (\lambda) =\lambda +i\delta,\,\lambda\in [0,+\infty)
\end{equation}
and $\delta$ 
a fixed value different from zero, represent the curve (\ref{LGAir}) with
\begin{equation}
r_0= \left (\frac32 |\delta| \right)^{2/3}
\end{equation} 
and with $\theta \in [0,\pi/3)$ if $\delta <0$ and  $\theta \in (\pi/3,2\pi/3]$ if $\delta >0$; similarly, $e^{-i\pi/3}z (\xi)$
gives a curve (\ref{LGAir}) with
$\theta\in (-\pi/3,0]$ if $\delta >0$ and $\theta\in [-2\pi/3,\pi/3]$ if $\delta<0$. Finally, $e^{\pm i\pi} z(\xi)$ 
gives the curves (\ref{LGAir}) for $\theta\in (\pi,2\pi/3)$ and  $\theta\in (-2\pi/3,-\pi)$. 

\subsection{{\normalsize Zeros of  {\boldmath$\cos\alpha {\rm Ai}(z)+\sin \alpha {\rm Bi}(z)$}}}

Denote 
\begin{equation}
\label{genai}
{\cal A}(\alpha,z)=\cos\alpha {\rm Ai}(z)+\sin\alpha {\rm Bi}(z).
\end{equation}
Using this combination and letting $\alpha$ be in general complex we can study the distribution
of zeros for any solution of the Airy equation. Indeed, given a solution $y(z)=\beta {\rm Ai}(z)+\gamma {\rm Bi}(z)$ the solution ${\cal A}(\arctan(\gamma/\beta),z)$ is proportional to $y(z)$ except when $\gamma/\beta=\pm i$; in this last case the zeros of $y(z)$ can be obtained
from those of ${\cal A}(\alpha,z)$ by letting $\Im \alpha \rightarrow \mp \infty$.

According to the description of the ASLs, we expect that zeros can appear which asymptotically tend
to one or several of the rays $\arg z =\pi, \pm  \pi/3$. Later, when we analyze the zeros of Bessel functions,
the zeros in the directions $\arg z =\pm  \pi/3$  will be particularly important. We analyze the distribution of zeros
in the three directions, starting with the negative real axis.

\subsubsection{Zeros approaching asymptotically the negative real axis}

For any real $\alpha$, the function ${\cal A}(\alpha,z)$ has an infinite number of negative real zeros, and
for any finite non real value of $\alpha$, there are zeros approaching the negative
real axis, but not on the axis. To see this, we write 
\begin{equation}
\label{asinAi}
\begin{array}{l}
{\rm Ai}(-z)\sim \pi^{-1/2}z^{-1/4}\left(\cos\varphi\, {\rm R}(\xi) +\sin\varphi\, {\rm S}(\xi)\right),\\
{\rm Bi}(-z)\sim \pi^{-1/2}z^{-1/4}\left(-\sin\varphi\, {\rm R}(\xi) +\cos\varphi\, {\rm S}(\xi)\right),
\end{array}
\end{equation}
with 
\begin{equation}
\label{varphi}
\varphi=\xi-\pi/4,\,\xi=\frac23 z^{3/2} 
\end{equation}
where ${\rm R}(\xi)$ and ${\rm S}(\xi)$ have the Poincar\'e-type expansions \cite[9.7.9, 9.7.11]{Olver:2010:ARF}
\begin{equation}
\label{asiPQ}
{\rm R}(\xi)=\displaystyle\sum_{k=0}^{\infty}(-1)^k \Frac{u_{2k}}{\xi^{2k}},\,
{\rm S}(\xi)=\displaystyle\sum_{k=0}^{\infty}(-1)^k \Frac{u_{2k+1}}{\xi^{2k+1}}
\end{equation}
in a sector $|\arg z|\le 2\pi/3-\delta$,
with $u_m$ numbers defined in \cite[9.7(i)]{Olver:2010:ARF}. The equality 
${\cal A}(\alpha ,-z)=0$ gives 
\begin{equation}
\label{firsas}
\Frac{{\rm Ai}(-z)}{{\rm Bi}(-z)}=-\tan\alpha .
\end{equation}
Using the leading terms in (\ref{asiPQ}) we have $\cot\varphi \sim \tan\alpha$ and inverting this relation we have
the estimation for the zeros
\begin{equation}
\label{zexi}
\xi\sim-\alpha-\Frac{\pi}{4}+k\pi,\, k\in {\mathbb Z},\,\Re \xi >0
\end{equation}
as asymptotic approximation for large positive $k$ ($\Re\xi<0$ does not make sense because (\ref{asiPQ})
is valid for $|\arg \xi|\le \pi-\epsilon$).

This corresponds to values of $-z$ approaching the negative real axis from above if $\Im \alpha >0$
and from below of $\Im \alpha <0$ as $k\rightarrow +\infty$; the approximation (\ref{zexi}) gives 
zeros which are on an ASL given by the LG approximation (see section \ref{ASLA}).

Eq. (\ref{zexi}) gives the dominant term in the asymptotic expansion for the negative real zeros of ${\cal A}(\alpha,z)$,
more terms can be obtained by inverting (\ref{firsas}) using additional terms of the expansions (\ref{asinAi}). For 
later use, it will be interesting to invert a more general relation than (\ref{firsas}), and we consider the asymptotic
inversion of
\begin{equation}
\label{firsas2}
\Frac{{\rm Ai}(-z)}{{\rm Bi}(-z)}=f(\alpha).
\end{equation}
If $f(\alpha)\neq \pm i$ and using (\ref{asinAi}) we obtain that 
\begin{equation}
\label{eqas}
\Frac{{\rm S}(\xi)}{{\rm R}(\xi)}= -\cot (\varphi - \chi ),\,\chi=\arctan (f(\alpha)).
\end{equation}

The case $f(\alpha)= \pm i$ would correspond to the zeros of ${\rm Ai}(-z) \pm i {\rm Bi}(-z)$, but considering 
\cite[9.2.11]{Olver:2010:ARF} we have ${\rm Ai}(-z) \pm i {\rm Bi}(-z)=2e^{\mp \pi i/3} {\rm Ai }(z e^{\pm i\pi /3})$, which does
not have zeros as $\Re z\rightarrow +\infty$ (all the zeros of ${\rm Ai}(z)$ are real and negative); we have in this case:
\begin{theorem}
\label{nozer}
Up to arbitrary multiplicative factors, ${\rm Ai }(z e^{\pm i\pi /3})$ are the only solutions of the Airy equation without
zeros as $\Re z \rightarrow -\infty$.
\end{theorem}

Inverting (\ref{eqas}) and using the expansions (\ref{asiPQ}) we have
\begin{equation}
\label{QP}
\xi+\Frac{\pi}{4}-\chi -k\pi=\arctan\Frac{{\rm S}(\xi)}{{\rm R}(\xi)}
\end{equation}
and following the steps \cite[section 7.6.1]{Gil:2007:NSF} we conclude that the zeros of (\ref{eqas}) as $\Re z\rightarrow
+\infty$ are given by
\begin{equation}
\label{zerAin}
z_k=T(t_k),\, t_k=\frac38 \pi (4k -1)+\frac32 \chi,\, k\in{\mathbb Z}\, \mbox{(}t_k>0\mbox{)}
\end{equation}
where $T(t)$ has an asymptotic expansion given in \cite[9.2.11]{Olver:2010:ARF} and additional terms can be found in \cite{Fab:1999:OTR}. The first terms of this expansion for large $t$ are
\begin{equation}
T(t)\sim t^{2/3}\left(1+\Frac{5}{48}t^{-2}-\Frac{5}{36}t^{-4}+\Frac{77125}{82944}t^{-6}-\cdots\right).
\end{equation}

Considering the particular case $f(\alpha)=-\tan\alpha$, we have that the zeros of ${\cal A}(\alpha,z)$ approaching the negative
real axis as $\Re z\rightarrow -\infty$ are given by
\begin{equation}
\label{firex}
z_k=-T(t_k),\, t_k=\frac38 \pi (4k -1)-\frac32 \alpha,\, k\in{\mathbb Z}\, \mbox{(}t_k>0\mbox{)}.
\end{equation}
Particular cases are the real zeros of ${\rm Ai}(z)$ ($\alpha=0$) and ${\rm Bi}(z)$ ($\alpha=\pi /2$), 
which correspond to Eqs. 9.9.6 and 9.9.10 of \cite{Olver:2010:ARF}. 

\begin{remark}
\label{remark1}
Observe that the smallest possible value
of $k$ in (\ref{firex}) may vary depending on the value of $\alpha$. For instance, 
comparing with Eqs. 9.9.6 and 9.9.10 of \cite{Olver:2010:ARF} it is clear that $k=1,2,...$ for $\alpha=0,\pi/2$; but
for $\alpha=-\pi/2$ (which algo gives the zeros of ${\rm Bi}(z)$) we should take $k=0,1,...$. 
 In order to be sure that the smallest possible positive value of $t_k$ is meaningful
we should make a detailed counting of zeros, using the argument principle as in \cite{Olver:1954:TAE}; we will not
consider this analysis here.
\end{remark}

The analysis of the zeros of the first derivative is very similar and leads to analogous results. One would start with Eqs.
9.7.10 and 9.7.12 of \cite{Olver:2010:ARF}, which we can write as
\begin{equation}
\begin{array}{l}
{\rm Ai}^{\prime}(-z)\sim \Frac{z^{1/4}}{\sqrt{\pi}}(\cos(\varphi-\pi/2)\tilde{R}(\xi)+\sin(\varphi-\pi/2)\tilde{S}(\xi)),\\
{\rm Bi}^{\prime}(-z)\sim \Frac{z^{1/4}}{\sqrt{\pi}}(-\sin(\varphi-\pi/2)\tilde{R}(\xi)+\cos(\varphi-\pi/2)\tilde{S}(\xi))
\end{array}
\end{equation}
where ${\rm R}(\xi)$ and ${\rm S}(\xi)$ have the Poincar\'e-type expansions 
\begin{equation}
\label{asiPQ2}
\tilde{R}(\xi)=\displaystyle\sum_{k=0}^{\infty}(-1)^k \Frac{v_{2k}}{\xi^{2k}},\,
\tilde{S}(\xi)=\displaystyle\sum_{k=0}^{\infty}(-1)^k \Frac{v_{2k+1}}{\xi^{2k+1}}
\end{equation}
in a sector $|\arg z|\le 2\pi/3-\delta$,
with $v_m$ numbers defined in \cite[9.7(i)]{Olver:2010:ARF}. Then, the equation
\begin{equation}
\Frac{{\rm Ai}^{\prime}(z)}{{\rm Bi}^{\prime}(z)}=f(\alpha)
\end{equation}
leads to 
\begin{equation}
\label{eqas2}
\Frac{\tilde{S}(\xi)}{\tilde{R}(\xi)}= -\cot (\varphi -\pi/2 - \chi ),\,\chi=\arctan (f(\alpha)).
\end{equation}
Therefore, the first approximation for the zeros of the derivatives can be obtained from those of the function by replacing
$\chi$ by $\chi +\pi/2$. In other words, in first approximation 
the location of the nodal curve with respect to the axis is the same for the zeros
and the zeros of the derivative, as it only depends on $\Im\chi$. This is consistent with the previous analysis based on the assumption that $A(z)$ is locally constant: the zeros of a given solution and those of the derivative lie
on the same ASL under this approximation.

For this reason, we will concentrate mainly on the zeros of the solutions of the ODEs, and not so much on the zeros
of the derivative. The results we will enunciate, except for specific asymptotic estimates, can be applied both for the function
and its derivative. For instance, Theorem \ref{nozer} is also true for the derivative and we have:
\begin{theorem}
\label{nozer2}
Up to arbitrary multiplicative factors, ${\rm Ai }(z e^{\pm i\pi /3})$ are the only solutions of the Airy equation without
zeros of its first derivative as $\Re z \rightarrow -\infty$
\end{theorem}
The same type of annotations can be made with respect to the results appearing in subsequent sections. 

\subsubsection{{\normalsize Zeros tending to the rays {\boldmath $\arg z=\pm \pi/3$} and general results}}
\label{argpi3}

If $\alpha\neq 0$, there is also an infinite number of zeros tending to the ray
 $\arg z=\pi/3$ (and similarly for $\arg z=-\pi/3$).
For $\alpha=0$ we are in the case of the Airy function ${\rm Ai}(z)$, which only has
negative real zeros.

Let us see, depending on $\alpha$, when there is a string of zeros above the ray
${\rm arg}\, z=\pi/3$ and when it lies below this ray.

Using \cite[9.2.11]{Olver:2010:ARF} we have
\begin{equation}
{\cal A}(\alpha,z)=e^{i\alpha}e^{-i\pi/3}{\rm Ai}(ze^{2\pi i/3})+e^{-i\alpha}e^{i\pi/3}{\rm Ai}(ze^{-2\pi i/3}),
\end{equation}
replacing $z$ by $z e^{i\pi/3}$ and using again \cite[9.2.11]{Olver:2010:ARF} (but with $z$ replaced by $z e^{-i\pi/3}$)
we obtain:
\begin{equation}
{\cal A}(\alpha,ze^{i\pi /3})=\left(e^{i\alpha}e^{-i\pi /3} +\frac12 e^{-i\alpha} e^{2i\pi /3}\right)\mbox{Ai}(-z)
-\Frac{i}{2}e^{-i\alpha}e^{2i\pi /3}\mbox{Bi}(-z).
\end{equation}
Because ${\rm Ai}(-z)$ and ${\rm Bi (-z)}$ do not have common zeros then ${\cal A}(\alpha,ze^{i\pi /3})$ and ${\rm Ai}(-z)$
do not have common zeros either and we have that if ${\cal A}(\alpha,ze^{i\pi /3})=0$ then
\begin{equation}
\label{aire2}
\Frac{{\rm Ai}(- z)}{{\rm Bi}(- z)}=\Frac{i}{1-2e^{2 i\alpha}}.
\end{equation}
Then if ${\cal A}(\alpha,ze^{i\pi /3})=0$ has zeros for $z$ real and positive (\ref{aire2}) holds, where the ratio
${\rm Ai}(- z)/{\rm Bi}(-z)$ is real; but 
this can only happen if $\Re (1-2 e^{2i\alpha})=0$ or, equivalently if $|1-e^{-2i\alpha}|=1$.

With respect to the zeros of ${\cal A}(\alpha,ze^{-i\pi /3})=0$, a similar analysis shows that this equation is equivalent to
\begin{equation}
\label{aire3}
\Frac{{\rm Ai}(- z)}{{\rm Bi}(- z)}=\Frac{-i}{1-2e^{-2 i\alpha}}
\end{equation}
and then there are positive real zeros of  ${\cal A}(\alpha,ze^{i\pi /3})=0$ only if $\Re (1-2 e^{-2i\alpha})=0$ or, 
equivalently, $|1-e^{2i\alpha}|=1$. On the other hand, ${\cal A}(\alpha,-z)=0$ has positive real zeros if and only if $\alpha$ is real. Therefore:
\begin{theorem}
${\cal A}(\alpha,z)$ has negative real zeros if and only if $\Im\alpha =0$ and it has 
zeros on the ray $\arg z=j \pi/3$ ($j=-1$ or $j=1$) if and only if $|1-e^{-2 j i\alpha}|=1$.

Any solution with zeros on a ray $\arg z=\pi$, $\arg z=\pi/3$ or $\arg z=-\pi/3$, is proportional to 
${\cal A}(\alpha,z)$, ${\cal A}(\alpha,z e^{-2i\pi/3})$ or ${\cal A}(\alpha,z e^{2i\pi/3})$ respectively for some real $\alpha$.
\end{theorem}

Also, as a consequence
\begin{theorem}
\label{sobreray}
Up to arbitrary multiplicative factors, ${\cal A}(\pm \pi/6,z)$ are the only solutions of the Airy equation with zeros
on the three principal anti-Stokes lines (the rays $\arg z=0,\pm \pi/3$).

\end{theorem}

On the other hand, similarly as in the case of Theorem \ref{nozer}, it is easy to see that when in (\ref{aire2}) we have
$1-2e^{2 i\alpha}=\pm 1$ ($\alpha=0$ or $\Im\alpha\rightarrow +\infty$), there are no zeros of ${\cal A}(\alpha,ze^{i\pi /3})$
as $\Re z \rightarrow +\infty$, and this corresponds to the fact that ${\rm Ai}(z)$ and ${\rm Ai}(z e^{-i\pi/3})$ are,
up to multiplicative factors, the only solutions without zeros approaching the ray $\arg z=\pi/3$. This is summarized
in the following theorem:

\begin{theorem}
All the solutions of the Airy equation have zeros tending asymptotically to the three rays $\arg z=-\pi,\pm\pi/3$ with the only
exception (up to multiplicative factors) of
${\rm Ai}(z)$, ${\rm Ai}(z e^{i\pi/3})$ and  ${\rm Ai}(z e^{-i\pi/3})$ which only have zero on one ray ($\arg z=\pi$, 
$\arg z=\pi/3$ and $\arg z=-\pi /3$ respectively).
\end{theorem}

The same arguments are true for the zeros of ${\cal A}'(\alpha,ze^{\pm i\pi /3})=0$ and the results of this section can be
easily adapted to these zeros. The starting point would be to replace ${\rm Ai}(- z)/{\rm Bi}(- z)$ by 
${\rm Ai}'(- z)/{\rm Bi}'(- z)$ in (\ref{aire2}) and (\ref{aire3}).

The only thing left to describe is the position of the zeros with respect to the principal ASLs when the zeros are not
on the rays themselves and to describe some asymptotic approximations. We already did this for the ray $\arg z=\pi$, 
and now we complete the analysis for the other two principal ASLs.

From (\ref{aire2}) we conclude that the zeros of ${\cal A}(\alpha,z e^{i\pi/3})$ as $\Re z\rightarrow +\infty$ are
given by (\ref{zerAin}) with
\begin{equation}
\chi=\arctan (i/(1-2 e^{2i\alpha}))=\Frac{i}{2}\log (1-e^{-2i\alpha}).
\end{equation}

Setting the right-hand side to zero in (\ref{zerAin}) gives the first approximation
\begin{equation}
\label{lalala}
\xi_k \sim -\frac12 \arg (1-e^{-2i\alpha})-\Frac{\pi}{4}+k\pi+\Frac{i}{2}\log\left|1-e^{-2i\alpha}\right|, k\in {\mathbb Z}\, 
\mbox{(}\Re\xi_k>0\mbox{)}.
\end{equation}
Therefore, there are zeros of ${\cal A}(\alpha,z e^{i\pi/3})$ above (below) the positive real axis if $|1-e^{-2i\alpha}|$ is
greater (smaller) than $1$. This means that the zeros of ${\cal A}(\alpha,z)$ approach the ASL $\arg z=\pi/3$ from above
if $|1-e^{-2i\alpha}|$ is greater than $1$ and from below if it is smaller.

The first approximation for the zeros of zeros of ${\cal A}(\alpha,z e^{-i\pi/3})$ as $\Re z\rightarrow +\infty$
can be obtained in a similar way. Now
\begin{equation}
\chi=\arctan (-i/(1-2 e^{-2i\alpha}))=\Frac{-i}{2}\log (1-e^{2i\alpha})
\end{equation}
which gives the first approximation
\begin{equation}
\label{lalala2}
\xi_k \sim \frac12 \arg (1-e^{2i\alpha})-\Frac{\pi}{4}+k\pi-\Frac{i}{2}\log\left|1-e^{2i\alpha}\right|, k\in {\mathbb Z}\, 
\mbox{(}\Re\xi_k>0\mbox{)}
\end{equation}
and the zeros approach the ASL $\arg z=-\pi/3$ from below (above)
if $|1-e^{-2i\alpha}|$ is greater (smaller) than $1$.

For the zeros of the derivative similar approximations hold, with $\chi$ replaced by $\chi+\pi/2$, and the location of the
zeros of the derivative with respect to the axis is the same.

From the discussion of section (\ref{ASLA}), we see that the values of $\tilde{z}= e^{i\pi/3} z =e^{i\pi/3}\left(\Frac{3}{2}\xi\right)^{2/3}$
corresponding to (\ref{lalala}) (that is, the approximations for the zeros of ${\cal A}(\alpha,z)$ approaching $\arg z=\pi/3$) 
lie exactly over the ASL (in the LG approximation) given by (\ref{LGAir}), with 
\begin{equation}
\label{erre0}
r_0=\left|\Frac{3}{4}\log|1- e^{-2i\alpha}|\right|^{2/3},
\end{equation}
where the angle $\theta$ in (\ref{LGAir}) is such that
 $\theta\in [0,\pi/3)$ if $|1-e^{-2i\alpha}|<1$ and $\theta\in (\pi/3,2\pi/3]$ if $|1-e^{-2i\alpha}|>1$.

Similarly, the first approximations for the zeros of ${\cal A}(\alpha,z)$ approaching $\arg z=-\pi/3$ lie exactly over 
the LG-ASL with 
\begin{equation}
\label{erre0}
r_0=\left|\Frac{3}{4}\log|1- e^{2i\alpha}|\right|^{2/3},
\end{equation}
and $\theta\in (-\pi /3,0]$ if $|1-e^{2i\alpha}|<1$ while $\theta\in (-2\pi/3,\pi/3]$ if $|1-e^{2i\alpha}|>1$.

The same is true for the zeros of the first derivative.

For the particular case of real $\alpha\in (0,\pi)$ we have the following: 
\begin{theorem}
The zeros of ${\cal A} (\alpha,z)=\cos\alpha {\rm Ai}(z)+\sin \alpha {\rm Bi}(z)$ are below (above) the ray 
$\arg z =\pi/3$ ($\arg z =-\pi /3$) if $\alpha\in (-\pi/6,\pi/6)$, above (below) 
this ray if  $\alpha\in [-\pi/2,-\pi/6)\cup (\pi /6,\pi/2]$ and 
exactly on the rays $\arg z=\pm \pi/3$ when $\alpha=\pm \pi/6$.
\end{theorem}
We notice that ${\cal A} (\pi /6,z e^{\pm i\pi /3})={\cal A} (\pi /6,-z)$ and 
${\cal A} (-\pi /6,z e^{\pm i\pi /3})=e^{\mp 2\pi i/3}{\cal A} (-\pi /6,-z)$, 
which shows explicitly that ${\cal A} (\pm \pi /6,z)$ 
have zeros over the rays $\mbox{arg}\, z=\pm \pi /3$ (see Theorem \ref{sobreray}). 

We stress again that Theorems 2.3 to 2.6 also hold for the first derivative.

The first approximation for the zeros is sufficient to have a clear picture of their distribution. However, 
proceeding as before for the zeros as $\Re z\rightarrow -\infty$,
we can obtain additional terms of the asymptotic expansions for large $k$. In particular,
${\cal A}(\alpha,z)$ has zeros approaching the ray $\arg z=\pi/3$ which are given by
\begin{equation}
z_k = e^{i\pi/3} T(t_k),\,t_k =\frac38 \pi (4k -1)+  i\frac34 \log(1- e^{-2 i\alpha}), k\in {\mathbb Z}\, 
\mbox{(}\Re t_k>0\mbox{)}.
\end{equation}
which for $\alpha=\pi/2$ gives the complex zeros of ${\rm Bi}(z)$ \cite[9.9.14]{Olver:2010:ARF}.

Similarly, ${\cal A}(\alpha,z)$ has zeros approaching the ray $\arg z=-\pi/3$ which are given by
\begin{equation}
z_k= e^{-i\pi/3} T(t_k),\,t_k =\frac38 \pi (4k -1)-  i\frac34 \log(1- e^{2 i\alpha}), k\in {\mathbb Z}\, 
\mbox{(}\Re t_k>0\mbox{)}.
\end{equation}

For the zeros of the derivative ${\cal A}'(\alpha,z)$ we have zeros at
\begin{equation}
\begin{array}{l}
z_k=-U(t_k),\, t_k=\frac38 \pi (4k -3)-\frac32 \alpha,\\
z_k = e^{i\pi/3} U(t_k),\,t_k =\frac38 \pi (4k -3)+  i\frac34 \log(1- e^{-2 i\alpha}),\\ 
z_k = e^{-i\pi/3} U(t_k),\,t_k =\frac38 \pi (4k -3)-  i\frac34 \log(1- e^{2 i\alpha})
\end{array}
\end{equation}
where $k\in{\mathbb Z}$ with $\Re t_k>0$ (see Remark \ref{remark1}). $U(t)$ has the asymptotic expansion 
\cite[9.9.19]{Olver:2010:ARF}; the first terms of the expansion are
\begin{equation}
U(t)\sim t^{2/3}\left(1-\Frac{7}{48}t^{-2}+\Frac{35}{288}t^{-4}-\Frac{181223}{207360}t^{-6}+\cdots\right).
\end{equation}

\section{{\large Zeros of Bessel functions}}

For the zeros of Bessel functions we follow a similar scheme as for Airy functions 
with one important difference. As before, we start by analyzing
the structure of the Stokes and anti-Stokes lines and then a more detailed analysis is performed using
asymptotics. But differently from the Airy case, asymptotics for large $z$ will be not enough to obtain a
complete picture of the distribution of zeros, and we will also need to consider uniform asymptotic approximations
for large order \cite{Olver:1954:TAE} in order to describe some of the zeros.

As happened for Airy functions, the picture of the zeros of the first derivative is very similar to that of the
zeros of the function, and their relative position with respect to the anti-Stokes lines the same. We are not
carrying a detailed description of the zeros of the first derivative, and we accept that the correspondence
can be made very easily.

\subsection{{\normalsize Stokes and anti-Stokes lines for the Bessel equation in the LG approximation}}
\label{StoBes}

We consider the differential equation satisfied by the Riccati-Bessel functions 
$\sqrt{z}J_{\nu}(z)$, with coefficient $A(z)=1-(\nu^2-1/4)/z^2$.  Let us focus on the case of Bessel functions 
of real orders $|\nu|>1/2$. The differential equation has two simple turning points at $z_{\pm}=\pm\sqrt{\nu^2-1/4}$ and because 
$\mbox{sign}(A'(z_{\pm}))=\mbox{sign}(z_{\pm})$ the principal 
ASLs emerge from $z_{\pm}$ at angles $0,\pm 2\pi/3$ for $z_{+}$ and $\pi ,\pm \pi/3$ for
$z_-$.

\begin{figure}[tb]
\vspace*{0.8cm}
\begin{center}
\centerline{\protect\hbox{\psfig{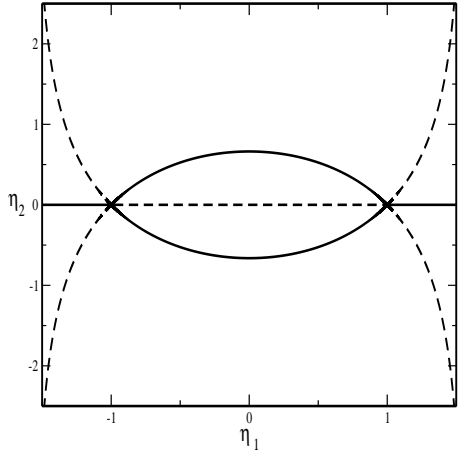}}}
\end{center}
\caption{The approximate principal anti-Stokes (dashed lines) and Stokes lines (solid lines) emerging from $\eta=\pm 1$ for the Bessel equation with $A(z)=1-(\nu^2 -1/4)/z^2$, $\nu$ real and $|\nu|>1/2$,  and $\eta =\eta_1 + i\eta_2 =z/\sqrt{\nu^2-1/4}$.
}
\label{fig1}
\end{figure}

 Writing the Riccati-Bessel equation $w''(z)+(1- \lambda^2/z^2)w(z)=0$, $\lambda=\sqrt{\nu^2 -1/4}$ in
the variable $\eta=z/\lambda$, the transformed
equation has coefficient $A(\eta)=\lambda^2 (1-1/\eta^2)$. 
For real $\lambda$ the ASLs emerging from the $\eta=\pm 1$ which are given by $\Im \int_{1}^z{\sqrt{A(u})}du=0$ can be written as
in \cite[Eq. (36)]{Segura:2012:CCZ} or equivalently as the curve in the complex plane $F=1$, where
\begin{equation}
\label{efe}
F(\eta)=\left|e^{\sqrt{1-\eta^2}}\Frac{\eta}{\sqrt{1-\eta^2}+1}\right|.
\end{equation}
The principal lines are shown in Fig. \ref{fig1}.  The eye-shaped region cuts the imaginary axis
 at $\eta=\pm i c$ where $c$ is the real root of 
$\sqrt{1+c^2}-\log((1+\sqrt{1+c^2})/c)=0$: 
\begin{equation}
\label{constant}
c= 0.66274321...
\end{equation}
The rest of ASLs are also given by $F=C$ with $C$ a constant; for $0<C<1$ we have closed curves inside the eye-shaped
region and for $C>1$ a curve for $\Im z>0$ with horizontal asymptotes as $\Re z\rightarrow \pm \infty$ and its complex
conjugated curve for $\Im z <0$.
$F<1$ gives the interior of the eye-shaped region and $F>1$ the exterior.

 According to the scheme of ASLs depicted in Fig. \ref{fig1} we expect that several types of zeros may appear 
(in the LG approximation): 

\begin{enumerate}
\item{}A string of zeros as $\Re z \rightarrow +\infty$ and with asymptote $\Im z =d$ for some real value $d$. Only 
one string (or none) of this type can be expected, either with positive or negative $d$.
\item{}String(s) of zeros as $\Re z \rightarrow -\infty$ and with asymptote $\Im z =d$ for some real value $d$;
two of such strings may appear, one with positive $d$ and one with negative $d$, but also one or none of these
strings is a possible situation. \footnote{Observe that the cases of zeros with $\Re z\rightarrow +\infty$ and 
$\Re z\rightarrow -\infty$ are different because the negative real axis is taken as the branch cut}
\item{}Zeros either inside or above/below the eye shaped region. If there are zeros inside the eye-shaped region,
they are located on a closed curve given by $F=C<1$ for some positive $C$. If the zeros
are outside ($F>1$), they may appear to lie on a same curve as the zeros as $\Re z\rightarrow \pm \infty$. Also, zeros
on the principal ASL forming the eye-shaped region are possible.
\end{enumerate} 

 This is coherent with the asymptotic behavior for the zeros of Bessel functions of $n$ integer 
(see \cite[10.21(ix)]{Olver:2010:BES}), but also for other cases described in the literature like Hankel functions
of real order \cite{Cruz:1982:ZHF} or the Bessel functions $J_{\nu}(z)$ and $Y_{\nu}(z)$ of semi-integer order 
\cite{Davies:1973:CZO}. For instance, the Hankel function $H^{(1)}_{\nu}(z)$ has zeros close to the eye-shaped region
for $\Im z <0$ (for an illustration for integer order see \cite[Fig. 10.21.4]{Olver:2010:BES}) 
and, depending on the value of $\nu$ (and as we will later see) an additional string of zeros below the cut
\cite{Cruz:1982:ZHF}. On the other hand the Bessel function $Y_n(z)$ of integer order 
has all the three types of zeros, with the zeros for $\Im z<0$ complex conjugated of those for $\Im z>0$; in fact, this is
true for $Y_{\nu}(z)$ with real orders and more generally also for $\cos\alpha J_{\nu}(z)-\sin\alpha Y_{\nu}(z)$ for some
values of $\alpha$, as we will later see.

In the following, we will provide information on the location of the zeros
for general solutions of the Bessel function of real orders, by considering the solutions
\begin{equation}
\label{general}
{\cal C}_{\nu}(\alpha,z)=\frac12 \left(e^{i\alpha} H^{(1)}_{\nu}(z)+ e^{-i\alpha} 
H^{(2)}_{\nu} (z)\right)=\cos\alpha J_{\nu}(z)-\sin\alpha Y_{\nu}(z).
\end{equation}
We can restrict the analysis to $\nu\ge 0$, because
using the connections formulas of \cite[10.4]{Olver:2010:BES}, one readily sees that
\begin{equation}
{\cal C}_{-\nu}(\alpha,z)={\cal C}_{\nu}(\alpha+\nu\pi,z).
\end{equation}
We will consider $\arg z\in (-\pi,\pi]$, that is, we will stay in the principal Riemann
sheet; the zeros in other Riemann sheets can be expressed in terms of zeros in the 
principal Riemann sheet using continuation formulas \cite[10.11]{Olver:2010:BES}.

We discuss in some detail the case of real $\alpha$, but letting $\alpha$ be complex
results for general combinations with complex coefficients will become available. It is
in particular convenient to analyze the limiting cases $\Im \alpha\rightarrow \pm \infty$,
which correspond to the zeros of Hankel functions $H^{(1)}(z)$ and $H^{(2)}(z)$.

\subsection{{\normalsize Zeros of Hankel functions}}

\label{Hankel}

Hankel functions $H^{(1)}(z)$ and $H^{(2)}(z)$ are the only
pair of independent solutions of the Bessel equation (up to constant
multiplicative factors) which 
do not have zeros for large and positive $\Re (z)$. The zeros of $H^{(2)}(z)$ are
complex conjugated of the zeros of $H^{(1)}(z)$.

That these functions have no zeros as $\Re (z)\rightarrow +\infty$ is obvious from
its asymptotic behaviour \cite[10.17.5-6]{Olver:2010:BES}. That the rest of solutions
have zeros as $\Re z\rightarrow +\infty$ follows
from this same asymptotic estimations, by considering the combination (\ref{general})
with general complex $\alpha$.
We have that, as $\Re z\rightarrow +\infty$ the equation ${\cal C}_{\nu}(\alpha,z)=0$ implies
that
\begin{equation}
e^{2i\omega}(1+{\cal O}(z^{-1}))\sim -e^{-2 i\alpha},\,\omega=z-\Frac{\nu\pi}{2}-\pi /4
\end{equation}
and inverting we have the estimation for the zeros
\begin{equation}
\label{largek}
z_k\sim -\alpha +\left(\Frac{\nu}{2}-\frac14+k\right)\pi,\, k\in {\mathbb Z}\,\mbox{(}\Re z_k>0\mbox{)}.
\end{equation}
Therefore, for any $\alpha\in {\mathbb C}$ we have zeros that run parallel to
the positive real axis with asymptote $\Im z =-\Im \alpha$; only if $\Im \alpha \rightarrow
\pm \infty$ these zeros do not exist (case of the Hankel functions).

Eq. (\ref{largek}) is the starting point for MacMahon asymptotic expansions \cite[10.21.19]{Olver:2010:BES};
in \cite[Example 7.9]{Gil:2007:NSF} additional terms in the expansion where given. The same expansions 
are valid in the general case with $\beta$ in  \cite[7.33]{Gil:2007:NSF} replaced by 
$\beta = -\alpha +\left(\Frac{\nu}{2}-\frac14+k\right)\pi$ (this parameter is denoted as $a$ in 
\cite{Olver:2010:BES}). 

Regarding the zeros close to the eye-shaped region, $H_{\nu }^{(1)}(z)$ has zeros on the
lower part ($\Im z<0$) but not on the upper part, and the contrary happens with
$H_{\nu}^{(2)}(z)$. These, together with $J_{\nu}(z)$, which only
has real zeros, are (up to multiplicative factors) the only solutions of the
Bessel equation that do not have zeros both with positive and negative imaginary part
in $|\Re z|<\nu$ for sufficiently large $\nu$. We postpone this analysis for section \ref{ojito}. 

Finally, $H^{(1)}_{\nu}(z)$ may have zeros below the negative real axis, depending on the value of $\nu$, while it does not
have zeros over the negative real axis.  

In the first place, is easy to see that there are no zeros over the negative real axis by 
using the continuation formula  \cite[10.11.5]{Olver:2010:BES}:
\begin{equation}
H_{\nu}^{(1)}(z e^{\pi i})=-e^{-\nu\pi i}H_{\nu}^{(2)} (z) .
\end{equation}
This implies that, because $H_{\nu}^{(2)} (z)$ does not have zeros as $\Re z\rightarrow +\infty$ with $\Im z<0$, then
$H_{\nu}^{(1)} (z)$ does not have zeros as $\Re z\rightarrow -\infty$ with $\Im z >0$.

With respect to the zeros below the negative real axis, 
considering the formula \cite[10.11.3]{Olver:2010:BES} 
we have
\begin{equation}
\label{con2}
H_{\nu}^{(1)}(z e^{-\pi i})=2\cos (\nu \pi)H_{\nu}^{(1)}(z)+e^{-\nu\pi i} H_{\nu}^{(2)}(z).
\end{equation}
Now from the Hankel's asymptotic expansions \cite[10.7.5-6]{Olver:2010:BES} for $H_{\nu}^{(1)}(z)$
and $H_{\nu}^{(2)}(z)$, valid for $|\arg\,z|\le \pi-\delta$, we conclude that the zeros if $H_{\nu}^{(1)}(z e^{-\pi i})=0$
for large $z$, $\Re z>0$, then
\begin{equation}
2\cos (\nu \pi)e^{i\omega}+e^{-\nu\pi i} e^{-iw}={\cal O}(z^{-1}),\omega=z-\Frac{\nu\pi}{2}-\Frac{\pi}{4},
\end{equation}
and, after inverting, we see that the zeros of $H_{\nu}^{(1)}(z e^{-\pi i})=0$ can be estimated as
\begin{equation}
z_k  \sim\Frac{\pi}{4}(1-2p)+k\pi +\Frac{i}{2}\log |2\cos\nu\pi|,k\in{\mathbb Z}\, \mbox{(}\Re z_s >0\mbox{)},
\end{equation}
where $p=0$ if if $\cos\nu\pi<0$ and $p=1$ if  $\cos\nu\pi>0$.

We observe that $\Im z_s>0$ if $|2\cos \nu \pi|>1$, which happens when $\{\nu\}\in [0,1/3)\cup (2/3,1)$, where 
$\{\nu\}=\nu-\lfloor \nu\rfloor$ is the fractional part of $\nu$. Only for these values does $H_{\nu}^{(1)}(z e^{-\pi i})$
have zeros as $\Re z\rightarrow +\infty$ with positive real part, and therefore only for these
values does $H_{\nu}^{(1)}(z )$ have an infinite string of zeros below the negative real axis. For $\{\nu\}\in (1/3,2/3)$,
 $H_{\nu}^{(1)}(z )$ no longer has zeros below the negative imaginary axis in the principal Riemann sheet because,
as $\{\nu\}$ becomes larger than $1/3$ (or smaller than $2/3$) these zeros are on the next Riemann sheet
($\arg(e^{-i\pi} z_k)<-\pi$). For an analysis of the trajectories of this type of zeros depending on the 
value of $\nu$ we refer to \cite{Cruz:1982:ZHF}.

In \cite{Cruz:1982:ZHF} it is also shown that the zeros below the negative real 
axis all cross the negative real axis precisely at $\{\nu\}=1/3,2/3$. A way to see this is by using (\ref{con2}) and writing 
the Hankel functions in terms of $J_{\nu}(z)$ and $Y_{\nu}(z)$; with this we see that $H_{\nu}^{(1)}(z e^{-\pi i})=0$ can
be written:
\begin{equation}
A J_{\nu}(z)+ B Y_{\nu}(z)=0,\,A=2\cos\nu\pi+e^{-\nu\pi i},\,B=i(2\cos\nu\pi-e^{-\nu\pi i}).
\end{equation}
Clearly $A\neq 0$, $B\neq 0$ and because $J_{\nu}(z)$ and $Y_{\nu}(z)$ are real for real positive $z$ and do not have common zeros, the previous equality for $z$ real and positive implies that $A/B\in {\mathbb R}$, which only occurs when $|\cos\nu\pi|=1/2$ and then 
$\{\nu\}=1/3,2/3$.

Therefore we conclude that only if $\{\nu\}\notin [1/3,2/3]$ there are zeros of $H^{(1)}(z)$ 
below the imaginary axis and with $\Im z_s  \rightarrow -\frac12 \log|2\cos\nu\pi|$ as $s\rightarrow +\infty$ ($\Re z_s\rightarrow
-\infty$). The same is true for the zeros of ${H^{(1)}}'(z)$.   

\subsection{{\normalsize Zeros of {\boldmath ${\cal C}_{\nu}(\alpha,z)=\cos\alpha J_{\nu}(z)-\sin \alpha Y_{\nu}(z)$}}}
%%5

Here we provide a description of the distribution of zeros for general cylinder functions 
${\cal C}_{\nu} (\alpha,z)=\cos\alpha J_{\nu}(z)-\sin\alpha Y_{\nu}(z)$, where $\alpha$ may be complex. We will analyze in more
detail some particular cases, like the case of real $\alpha$, but the general results
provide a complete picture on the distribution of zeros. 

\subsubsection{{\normalsize Zeros parallel to the positive real axis as $\Re z\rightarrow +\infty$ and MacMahon-type expansions}}

As explained before, for any $\alpha$ there is a string of zeros running parallel to the real axis as $\Re z\rightarrow +\infty$
which are approximately given by (\ref{largek}).
Only as $\Im\alpha \rightarrow \pm \infty$ these zeros do not exist (case of the Hankel functions). We briefly summarize how
asymptotic expansions can be obtained for these zeros; the same procedure can be applied to the zeros running parallel
to the negative real axis (if any), that we will consider in the next section.

Inverting ${\cal C}_{\nu}(\alpha,z)=0$ is equivalent to
inverting $Y_{\nu}(z)/J_{\nu}(z)=\cot\alpha$; we write $\chi =\pi/2-\alpha$ and then we consider
the inversion of 
\begin{equation}
\label{nolnolnol}
\Frac{Y_{\nu}(z)}{J_{\nu}(z)}=\tan\chi
\end{equation}
as $\Re z\rightarrow +\infty$. We can let $\chi$ be complex in general and therefore obtain expansions for general
combinations of Bessel functions (with Hankel functions as the limits 
$\Im \chi\rightarrow \pm \infty$). Now, we write \cite[10.4]{Olver:2010:BES}
\begin{equation}
\label{PoinB}
\begin{array}{l}
J_{\nu}(z)=\left(\Frac{2}{\pi z}\right)^{1/2}(\cos\omega {\rm P}(\nu,z) -\sin\omega {\rm Q}(\nu,z)),\\
Y_{\nu}(z)=\left(\Frac{2}{\pi z}\right)^{1/2}(\sin\omega {\rm P}(\nu,z) +\cos\omega {\rm Q}(\nu,z)),
\end{array}
\end{equation}
where $\omega=z-\Frac{\nu\pi}{2}-\Frac{\pi}{4}$.
${\rm P}(\nu,z)$ and ${\rm Q}(\nu,z)$ have asymptotic expansions in inverse powers of $z$ for $|\arg z|<\pi-\delta$; 
 ${\rm P}(\nu,z)={\cal O}(1)$ while ${\rm Q}(\nu,z)={\cal O}(z^{-1})$. 

With this, the equation (\ref{nolnolnol}) becomes
\begin{equation}
\label{qupe}
\Frac{{\rm Q}}{{\rm P}}=\tan (\chi-\omega),
\end{equation}
and inverting this relation gives
\begin{equation}
z_s=\left(s+\nu/2-\Frac{3}{4}\right)\pi+\chi-\arctan\left(\Frac{{\rm Q}}{{\rm P}}\right),s\in {\mathbb Z}\, \mbox{(}
\Re z_s >0 \mbox{)}.
\end{equation}
Neglecting the last term ($Q/P={\cal O}(z^{-1})$) we have our first approximation, and by re-substitution we can generate as
many terms as needed of the asymptotic expansion for large $s$ 
(see \cite[Example 7.9]{Gil:2007:NSF}). This gives the MacMahon expansions 
\begin{equation}
\label{MMO}
z_s\sim \beta -\Frac{\mu -1}{2}\displaystyle\sum_{i=0}^{\infty}\Frac{p_i (\mu)}{(4\beta)^{2i+1}}
\end{equation}
where $p_{i}(\mu )$ are certain polynomials of $\mu$ of degree $i$. The parameter $\beta$ is given by
\begin{equation}
\label{beta}
\beta=\left(s+\Frac{\nu}{2}-\Frac{3}{4}\right)\pi+\chi .
\end{equation}
In particular, if we are considering the zeros of ${\cal C}_{\nu}(\alpha,z)=\cos\alpha J_{\nu}(z)-\sin\alpha Y_{\nu}(z)$
as $\Re z\rightarrow \infty$, then $\chi =\pi/2-\alpha$ and
\begin{equation}
\label{betaeta}
\beta=\left(s+\Frac{\nu}{2}-\Frac{1}{4}\right)\pi-\alpha.
\end{equation}
And for $\alpha=0$ and $\alpha=\pi/2$ we are in the cases of the positive zeros of $J_{\nu}(z)$ and $Y_{\nu}(z)$, and the
expansions correspond to \cite[10.21.19]{Olver:2010:BES} (where $\beta$ is denoted as $a$). For any other real $\alpha$,
we have positive real zeros, while for complex $\alpha$ the zeros as $\Re z\rightarrow +\infty$ have asymptote $\Im z =-\Im \alpha$. 

The procedure for inverting (\ref{nolnolnol}) can be used to
 provide MacMahon-type expansions also for the case of the zeros parallel to the negative real axis (if any) 
 just by replacing $\chi$ by its corresponding value.

For the zeros of the derivatives the analysis is very similar. We would start with \cite[10.17.9-10]{Olver:2010:BES}, which
we can write as
\begin{equation}
\label{PoinB}
\begin{array}{l}
J_{\nu}'(z)=\left(\Frac{2}{\pi z}\right)^{1/2}(\cos(\omega +\pi/2) \tilde{P}(\nu,z) 
-\sin(\omega +\pi/2) \tilde{Q}(\nu,z)),\\
Y_{\nu}'(z)=\left(\Frac{2}{\pi z}\right)^{1/2}(\sin(\omega +\pi/2) \tilde{P}(\nu,z) +
\cos(\omega +\pi/2) \tilde{Q}(\nu,z)),
\end{array}
\end{equation}
where, as before $\omega=z-\Frac{\nu\pi}{2}-\Frac{\pi}{4}$.
$\tilde{P}(\nu,z)$ and $\tilde{Q}(\nu,z)$ have asymptotic expansions in inverse powers of $z$ for  $|\arg z|<\pi-\delta$; 
 $\tilde{P}(\nu,z)={\cal O}(1)$ while $\tilde{Q}(\nu,z)={\cal O}(z^{-1})$. 

With this, the equation (\ref{nolnolnol}) becomes
\begin{equation}
\Frac{\tilde{Q}}{\tilde{P}}=\tan (\chi-\pi/2-\omega).
\end{equation}
The situation is analogous to (\ref{qupe}) but with $\chi$ replaced by $\chi-\pi/2$. It is then simple to write
down the corresponding expansion by using the MacMahon expansions \cite[10.21.20]{Olver:2010:BES}

Because it is the imaginary part of the zeros which controls the relative position with respect to the real axis, 
the results for the zeros of ${\cal C}(\alpha,z)$ and its first derivative are the same in
this sense (because $\chi$ is shifted
by a real amount); we will not consider again explicitly the case of the first derivative.

\subsubsection{{\normalsize Zeros parallel to the branch cut as {\boldmath $\Re z\rightarrow -\infty$}}}
\label{parbra}

Similarly as happened for the Hankel functions, for estimating the zeros running parallel to the branch cut as $\Re z\rightarrow -\infty$, we need to use continuation formulas so that the expansions (\ref{PoinB}) are applicable.

Using 
\cite[11.11.1-2]{Olver:2010:BES} we have 
\begin{equation}
{\cal C}_{\nu}(\alpha,z e^{im\pi})=(\cos\alpha e^{im\nu\pi}-2i\sin\alpha\sin(m\nu \pi)\cot (\nu\pi))J_{\nu}(z)-\sin\alpha e^{-im\nu\pi}
Y_{\nu}(z).
\end{equation}
We will take $m=\pm 1$.

Then, ${\cal C}_{\nu}(\alpha,z e^{im\pi})=0$ implies Eq. (\ref{nolnolnol}) with
\begin{equation}
\label{ratio}
\tan\chi =\cot\alpha e^{2 i m\nu\pi}-2 i e^{i m\nu \pi}\sin(m\nu\pi) \cot (\nu\pi) =A+iB,
\end{equation}
where
\begin{equation}
\begin{array}{l}
A=\cot\alpha +2 \sin^2 (m\nu\pi)(\cot(\nu\pi)-\cot\alpha),\\ 
B=-sin(2 m\nu\pi)(\cot(\nu\pi)-\cot\alpha).
\end{array}
\end{equation}

Taking only the first term in Eq. (\ref{MMO}) we have as an estimation for the values
of $z$ satisfying ${\cal C}_{\nu}(z e^{im\pi})=0$:
\begin{equation}
\label{zsf}
z_k\sim \left(k+\Frac{\nu}{2}-\Frac{3}{4}\right)\pi+\arctan(A+iB), k\in {\mathbb Z}\, 
\mbox{(}\Re z_k>0\mbox{)}.
\end{equation}
The complete expansion would be (\ref{MMO}) where now $\chi=\arctan(A+iB)$.

Using the first approximation (\ref{zsf}) we have 
\begin{equation}
\Im z_k \sim \Frac{1}{2}\log\left|\Frac{A+(B+1)i}{A+(B-1)i}\right|=a
\end{equation}
and there is an infinite number of zeros of ${\cal C}_{\nu}(\alpha,z e^{im\pi})=0$
for $\Re  z >0$ with asymptote $\Im z=a$.

Now, if we take $m=1$ and it turns out that $a<0$ then there is an infinite number of zeros of ${\cal C}_{\nu}(\alpha,z e^{im\pi})$ 
with $z$ below the positive real axis; therefore, ${\cal C}_{\nu}(z)$
will have an infinite number of zeros above the negative real axis. Contrary, if $a>0$ 
 the zeros of ${\cal C}_{\alpha,\nu}(z)$ would be on the next Riemann sheet after turning
with angle $\pi$ (because then $\arg (e^{i\pi}z_k)>\pi$).

Let us consider the particular case of $\alpha$ real, $\alpha\in [0,\pi)$. Observe that for real parameters 
$\mbox{sign}(a)=\mbox{sign} (B)$. 
Therefore, we have zeros over the branch cut (and also the conjugated zeros below the branch cut) only if
$B<0$. Then the set of values for which there exist zeros above (and below)
the branch cut is:
\begin{enumerate}
\item{}$\alpha\neq 0$ if $\nu\in {\mathbb Z}$
\item{}$\alpha >\{\nu\}\pi$ if $\{\nu\}\in (0,1/2)$
\item{}$\alpha <\{\nu\}\pi$ if $\{\nu\}\in (1/2,1)$
\end{enumerate}

We note that as, $\alpha$ varies,
the zeros of ${\cal C}_{\nu}(z)$ can be at the branch cut only if $B=0$, because ${\cal C}(\alpha,e^{im\pi}z)=0$ implies that
$$
\Frac{Y_{\nu}(z)}{J_{\nu}(z)}=A+Bi
$$ 
and $\Frac{Y_{\nu}(z)}{J_{\nu}(z)}$ is real for real and positive $z$. Therefore, if the zeros parallel to the branch cut
as $\Re z\rightarrow -\infty$ cross the negative real axis as $\alpha$ varies, they will all cross this axis for a same value of 
$\alpha$. The zeros lie over the negative real axis in the following cases:
\begin{enumerate}
\item{}$\alpha=0$
\item{}$\{\nu\}$ half-odd
\item{}$\{\nu\}$ not a half-odd number and $\alpha = \{\nu\}\pi$
\end{enumerate}
where $\{\nu\}$ is the fractional part of $\nu$.

For the particular case of integer order, there exist zeros over and below the branch cut for any $\alpha\neq 0$,
with asymptotes as $\Re z\rightarrow -\infty$ given by
\begin{equation}
\Im z =\pm \Frac{1}{4}\log(1+8\sin^2\alpha).
\end{equation}
And for $\alpha=0$, the result 10.21.46 of \cite{Olver:2010:BES} is reproduced.

For the more general case of complex parameters, there is a string of zeros of ${\cal C}_{\nu}(\alpha,z)$ 
running above
the negative real with asymptote
\begin{equation}
\label{paraa}
\Im z =-\Frac{1}{2}\log\left|\Frac{A+(B+1)i}{A+(B-1)i}\right|=-a
\end{equation}
if $a<0$. Similarly (taking $m=-1$ in the previous analysis), there is a string of zeros running below
the negative real with asymptote
\begin{equation}
\label{parab}
\Im z =-\Frac{1}{2}\log\left|\Frac{A+(-B+1)i}{A-(B+1)i}\right|=-b
\end{equation}
if $b>0$.

The same is true for the zeros of the first derivative.

\subsubsection{{\normalsize Zeros close to the eye-shaped region or inside the eye-shaped region}}
\label{ojito}

 The zeros close to the eye-shaped region are not so simple to analyze. Poincar\'e asymptotics
is of no use here, but uniform asymptotics for large $\nu$ gives a clear picture. 
In 1954, Olver developed powerful uniform asymptotic expansions for Bessel functions of large
order \cite{Olver:1954:TAE}, uniformly valid for $\tilde{z}=z/\nu \in (0,+\infty)$ as the order goes to infinity. These are expansions with Airy
functions as main approximants, and they read \cite[10.20.4-5]{Olver:2010:BES}:
\begin{equation}
\begin{array}{l}
J_{\nu}(\nu \tilde{z})\sim \left(\Frac{4\zeta}{1-\tilde{z}^2}\right)^{1/4}\nu^{-1/3}
\left[{\rm Ai}\left(\nu^{2/3}\zeta\right)\displaystyle\sum_{k=0}^{\infty}\Frac{A_{k}(\zeta)}{\nu^{2k}}
+{\rm Ai}^{\prime}\left(\nu^{2/3}\zeta\right)\displaystyle\sum_{k=0}^{\infty}\Frac{B_{k}(\zeta)}{\nu^{2k+1}}
\right],\\
Y_{\nu}(\nu \tilde{z})\sim -\left(\Frac{4\zeta}{1-\tilde{z}^2}\right)^{1/4}\nu^{-1/3}
\left[{\rm Bi}\left(\nu^{2/3}\zeta\right)\displaystyle\sum_{k=0}^{\infty}\Frac{A_{k}(\zeta)}{\nu^{2k}}
+{\rm Bi}^{\prime}\left(\nu^{2/3}\zeta\right)\displaystyle\sum_{k=0}^{\infty}\Frac{B_{k}(\zeta)}{\nu^{2k+1}}
\right]
\end{array}
\end{equation}
with $\zeta$ the solution of the differential equation
\begin{equation}
\label{defzeta}
\left(\Frac{d\zeta}{d\tilde{z}}\right)^2 =\Frac{1-\tilde{z}^2}{\zeta \tilde{z}^2}
\end{equation}
that is infinitely differentiable on the interval $0<\tilde{z}<\infty$, including $\tilde{z}=1$, and continued analytically to the 
complex plane cut along the negative real axis. 

From these approximations, we observe that for large enough $\nu$ the zeros of Bessel functions 
${\cal C}_{\nu}(\alpha,z)=\cos\alpha J_{\nu}(z)-\sin\alpha Y_{\nu}(z)$, $z=\nu\tilde{z}$, are given, in first
approximation, by the zeros of ${\cal A}(\alpha,\nu^{2/3}\zeta)=\cos\alpha {\rm Ai}(\nu^{2/3}\zeta)
+\sin\alpha {\rm Bi}(\nu^{2/3}\zeta)$. Therefore, for
$\nu$ sufficiently large, we can infer the distribution of the zeros of Bessel functions from that of Airy functions that we studied before.

Similarly, we can deduce the distribution of zeros of derivative ${\cal C}'_{\nu}(\alpha ,z)$ from
that of ${\cal A}'(\alpha,\nu^{2/3}\zeta)$ (see [10.20.7-8]\cite{Olver:2010:BES}). And because the zeros of 
Airy functions and of their first derivative  lie on the same curves (in first approximation),
the discussion for the zeros of ${\cal C}'_{\nu}(\alpha ,z)$ will be similar to
that of ${\cal C}_{\nu}(\alpha ,z)$. 

In this section we will not be concerned with the analysis of the asymptotic expansions of the
zeros for large $\nu$. Our main goal is to determine when the zeros related to the eye-shaped region exist and to 
describe the curve (ASL) where they lie. In particular, we will determine the intersection of this ASLs with the 
imaginary axis, and
also with the real axis for curves inside the eye-shaped region. As we discuss later, this is an important piece of 
information for the numerical computation of these zeros, particularly the intersection with the imaginary axis.

\begin{figure}[tb]
\vspace*{0.8cm}
\begin{center}
\begin{minipage}{3cm}
\centerline{\protect\hbox{\psfig{file=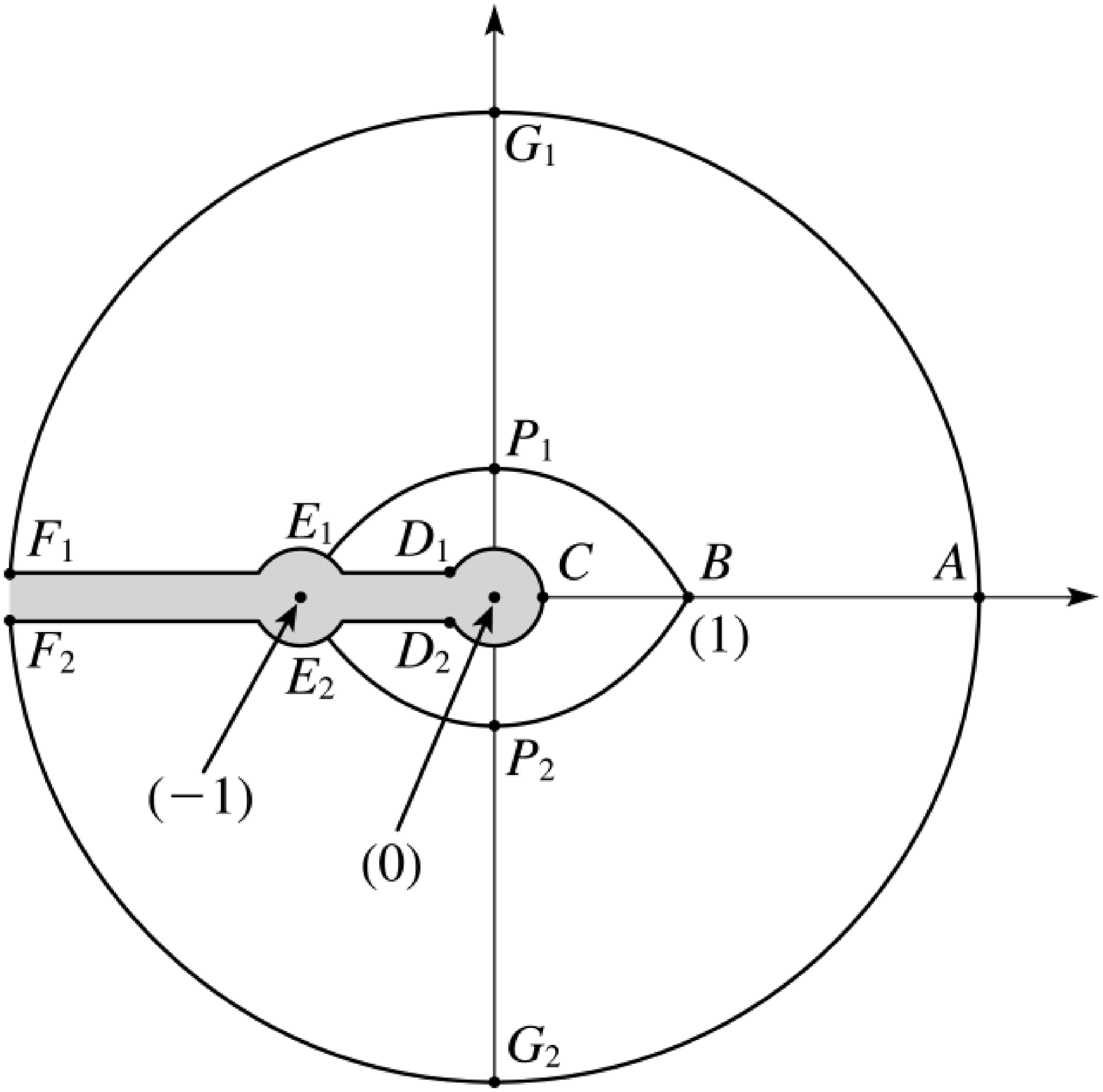,angle=0,height=5.5cm,width=5.5cm}}}
\end{minipage}
\hspace*{3cm}
\begin{minipage}{3cm}
\centerline{\protect\hbox{\psfig{file=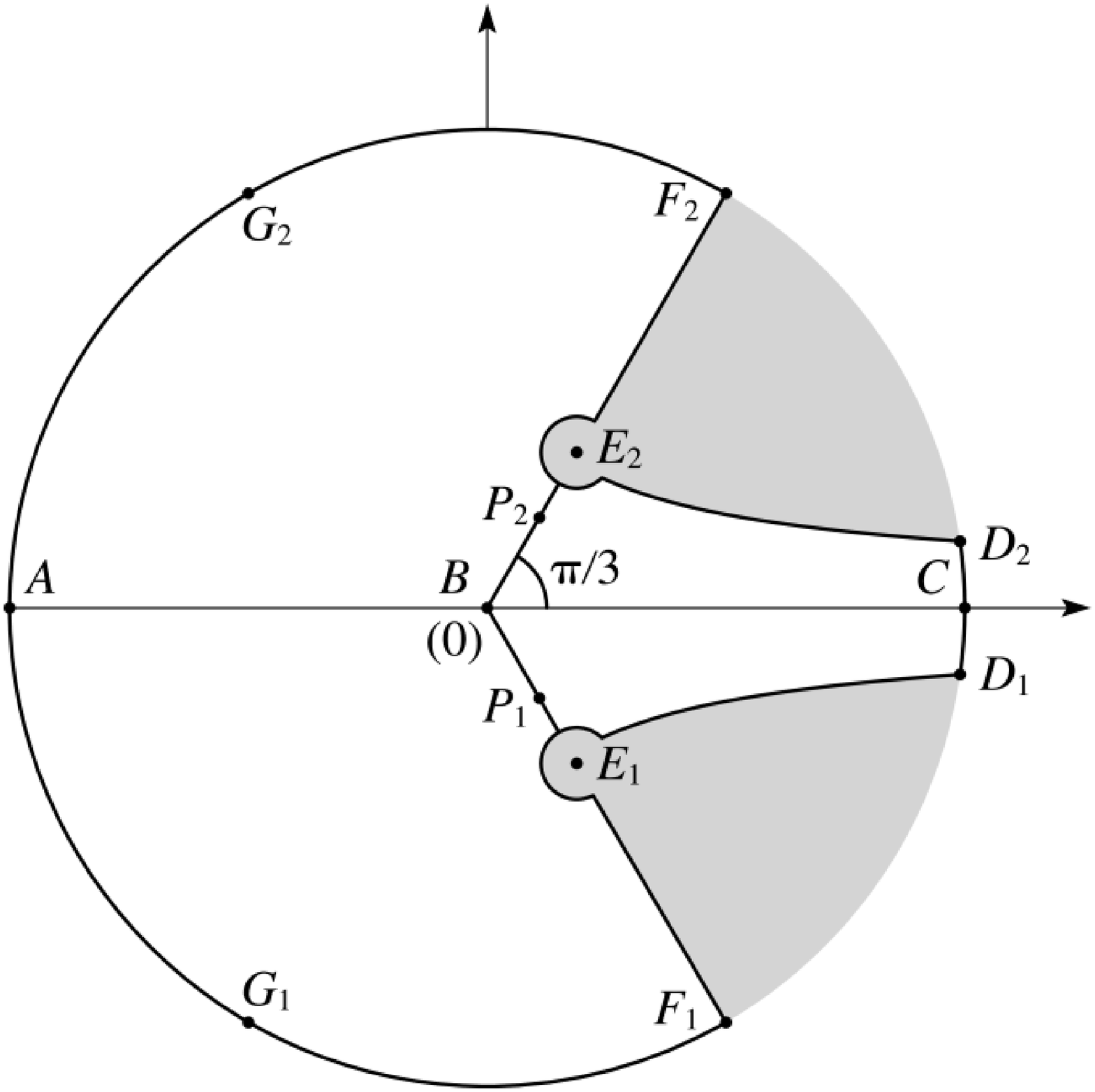,angle=0,height=5.5cm,width=5.5cm}}}
\end{minipage}
\end{center}
\caption{$\tilde{z}$-domain (left) and $\zeta$-domain (right), with corresponding points. These figures were taken from 
http://dlmf.nist.gov/10.20; they are copyrighted by NIST and used with permission.}
\label{fig2}
\end{figure}

 It is important to bear in mind the correspondence between the values of $\tilde{z}$ and those of $\zeta$ (Figure
\ref{fig2}). Observe that the curves $BP_1E_1$ and $BP_2E_2$ in the $\tilde{z}$-plane (the limits of the eye shaped region)
correspond to the line segments 
\begin{equation}
\zeta = e^{\pm i\pi /3}\tau\,,0\le\tau\le \left(\frac32 \pi\right)^{2/3}
\end{equation}
respectively. Therefore, for analyzing the zeros close to the eye-shaped region we need to 
consider the zeros of Airy functions as $\arg z\rightarrow \pm \pi/3$, $|z|\rightarrow \infty$ 
(section \ref{argpi3}). We consider the zeros approaching $\arg z =\pi /3$, which give the zeros corresponding
to the lower part of the eye-shaped region; the description for the upper part is analogous but considering the
direction $\arg z=-\pi/3$. 

We note the similarity between the domains determined by the ASLs in section \ref{StoBes} and the different regions in the $\tilde{z}$-domain of Fig. \ref{fig2}. The only difference is that the eye-shaped curve defined from (\ref{efe}) and setting 
$F(\eta)=1$ is given
in terms of $\eta=z/\sqrt{\nu^2-1/2}$, while the equation for the curve $E_1 P_1 B P_2 E_2$ of Fig.2 is given in terms of 
$\tilde z=z/\nu$ ($F(\tilde{z})=1$). In practice, those curves are very similar if $\nu$ is not small, and we expect that the
information we will obtain from asymptotics is consistent with the LG approximation. The eye-shaped domain inside
the curve $E_1 P_1 B P_2 E_2$, which  will be denoted as $K$, is 
\begin{equation}
K=\{\tilde{z}\in {\mathbb C}: F(\tilde{z})<1\}
\end{equation}
with $F$ given 
by (\ref{efe}).

In the $\zeta$ plane, we will have zeros below (above) the ray $\arg \zeta=\pi/3$ if $|1-e^{-2i\alpha}|$ is smaller (greater)
than $1$ (section \ref{argpi3}). 
If the zeros are below the ray $\arg \zeta=\pi/3$, this gives a finite number of zeros because the boundary 
given by the curve $D_2 E_2$ in the $\zeta$-domain would be reached by the string of Airy zeros as $|z|$ becomes large; this
corresponds to zeros inside the eye-shaped region and is in agreement with the picture given by the LG approximation. If
the zeros are above the ray $\arg z=\pi/3$ then, in the LG approximation for Airy functions, 
we have an ASL with an infinite number of zeros
both as $\arg \zeta\rightarrow \pi/3$ and as $\arg \zeta \rightarrow \pi$; as $\arg \zeta \rightarrow \pi$, in the $z$ variable
this would give a string of zeros asymptotically parallel to the positive real axis while, as $\arg \zeta\rightarrow \pi/3$,
we have a string of zeros below the negative real axis. This however, is an approximate picture and, as we have already
shown, the string of zeros as $\arg \zeta\rightarrow \pi/3$ may be absent in the principal Riemann sheet; it is important
to bear in mind that in the $\zeta$-domain there is a branch cut given by the ray $\arg z=\pi/3$ from the point $E_2$ to infinity (there is also the cut $\arg z=-\pi/3$ from the point $E_1$ to infinity).

In any case, there will be always zeros with $\Re \tilde{z}\in (-1,1)$  and $\Im \tilde{z}<0$ if $\nu$ is sufficiently large 
(in fact even for $\nu$ quite small)
whenever $\log  |1-e^{-2i\alpha}|$ is finite; therefore, for large enough $\nu$ there are always zeros satisfying these
conditions except when $|1-e^{-2i\alpha}|=0,\infty$, which is the case of the first kind Bessel
function (for $\alpha=0$ ${\cal C}_{\nu}(0,z)=J_{\nu}(z)$) and the Hankel function $H^{(2)}(z)$ 
(as $\Im\alpha\rightarrow +\infty$, $\tan\alpha\rightarrow -i$ corresponding to $H^{(2)}(z)=J_{\nu}(z)-iY_{\nu}(z)$). Similar
arguments can be used for the zeros with  $\Re \tilde{z}\in (-1,1)$  and $\Im \tilde{z}>0$. We summarize this in the 
following result.

\begin{theorem}
Let $\nu$ be real and positive and sufficiently large. 
The functions $J_{\nu}(z)$, $H_{\nu}^{(1)}(z)$ and $H_{\nu}^{(2)}(z)$ are, up to constant factors,
the only three solutions of the Bessel equation which do not have zeros for $\Re \tilde{z} \in (-1,1)$, 
$\tilde{z}=z/\nu$, both for $\Im z>0$ and $\Im z<0$
simultaneously.
 
$J_{\nu}(z)$ has no zeros such that $\Re z\in(-\nu,\nu)$.

$H_{\nu}^{(1)}(z)$ has zeros such that $\Re z\in(-\nu,\nu)$ and $\Im z<0$.

$H_{\nu}^{(2)}(z)$ has zeros such that $\Re z\in(-\nu,\nu)$ and $\Im z>0$.

\end{theorem}

For brevity, when a Bessel function ${\cal C}_{\nu}(\alpha,\nu \tilde{z})$ has zeros for $\Re \tilde{z} \in (-1,1)$ 
we will say that it has Airy-type zeros.

If $|1-e^{-2i\alpha}|<1$ we have zeros inside the domain $K$ for $\Im \tilde{z}<0$ and the ASL containing these
zeros intersects both the real axis and the imaginary axis. Contrarily, if $|1-e^{-2i\alpha}|>1$ the ASL corresponding
to these zeros will only intersect the imaginary axis. 

We start with the case of zeros inside $K$ and with negative real part, that is $|1-e^{-2i\alpha}|<1$.
From previous discussion, we know that these
zeros lie on the curve
\begin{equation}
\nu^{2/3}\zeta = r_0 \left| \cos \Frac{3\theta}{2}\right|^{-2/3} e^{i\theta},\, r_0=\left|\frac34 \log |1-e^{-2i\alpha}| \right|^{2/3},\theta\in (0 ,\pi /3)
\end{equation}
under the LG approximation.

Now, considering that \cite[10.20.2]{Olver:2010:BES}
\begin{equation}
\zeta=\left(\Frac{3}{2}f(\tilde{z})\right)^{2/3},f(\tilde{z})=\log\left(\Frac{1+\sqrt{1-\tilde{z}^2}}{\tilde{z}}\right)-
\sqrt{1-\tilde{z}^2}
\end{equation}
we can put in correspondence $\theta$ and $\tilde{z}$ as follows,
\begin{equation}
\label{mastereye}
\nu f(\tilde{z}) =\frac12 |\log |1- e^{-2 i\alpha}||\left(1+i\tan\Frac{3\theta}{2}\right),
\end{equation}
where we have used that $\cos \Frac{3\theta}{2}>0$.

When $\theta=0$, $\zeta$ becomes real and in $(0,1)$ and the corresponding value of $\tilde{z}$ gives the cut of the curve
containing the zeros inside the eye-shaped region with the real axis.
We have that the curve of zeros of ${\cal C}_{\nu} (\alpha, z)$ cuts the positive real $x$-axis 
at $\nu \tilde{x}_0$ with $\tilde{x}_0$
the positive real root of
\begin{equation}
\label{xeq}
\begin{array}{l}
f(\tilde{x}_0)=-\Frac{1}{2\nu}\log|1-e^{-2i\alpha}|,\\
\\
f(\tilde{x}_0)=\log\left(\Frac{1+\sqrt{1-\tilde{x}_0^2}}{\tilde{x}_0}\right)-\sqrt{1-\tilde{x}_0^2}
\end{array}
\end{equation}
and it also cuts the negative real axis at $-\nu \tilde{x}_0$.

As $\theta$ increases from $\theta=0$, the curve in the $\tilde{z}$ plane moves towards the negative imaginary real axis until a value
of $\theta$ in $(0,\pi/3)$ is reached such that becomes purely imaginary. Setting $\tilde{z}=-i \tilde{y}_0$ 
in (\ref{mastereye}) the real
part gives:

\begin{equation}
\label{yeq}
g(\tilde{y}_0)=-\Frac{1}{2\nu}\log|1-e^{-2i\alpha}|,\,g(\tilde{y}_0)=\log\left(\Frac{1+\sqrt{1+\tilde{y}_0^2}}
{\tilde{y}_0}\right)-\sqrt{1+\tilde{y}_0^2}.
\end{equation}

In a similar way, it is easy to check that in the case $|1-e^{-2i\alpha}|>1$, the cut with the imaginary axis $\tilde{z}=
-i\tilde{y}_0$ is 
also given by the solution of (\ref{yeq}).
For the case $|1-e^{-2i\alpha}|=1$ the solution of this equation ($g(\tilde{y}_0)=0$) gives the cut of the boundary 
$F(\tilde{z})=1$ with the imaginary axis
($-i y_0$).
For the Airy-type zeros with positive imaginary part (corresponding to $\arg \zeta \rightarrow -\pi/3$)
the same equations but replacing $\alpha$ by $-\alpha$ give the cut with the axis. 

Summarizing, the cuts of the curve containing
the Airy-type zeros of ${\cal C}_{\nu}(\alpha, z)$ with the
imaginary axis are $j i \nu \tilde{y}_j$, $j=\pm 1$, with $\tilde{y}_j$ the (positive) solution of 
\begin{equation}
\label{masterojito}
g(\tilde{y}_j)+\Frac{1}{2\nu}\log|1-e^{2 j i\alpha}|=0.
\end{equation} 
We observe that the function $g(y)$ is monotonic and $g(0^+)=+\infty$, $g(+\infty)=-\infty$, therefore the value of the
cut with the positive (or the negative) imaginary axis, as expected, can be any positive real value because 
$\Frac{1}{2\nu}\log|1-e^{\pm 2i\alpha}|$ can take any possible real values for complex $\alpha$.
However, as $\nu\rightarrow +\infty$ for a fixed $\alpha$ all the zeros tend to cluster on the 
boundary of the domain $K$.

The approximation of the point of intersection of the curve containing the Airy-type zeros with the imaginary axis 
given by (\ref{masterojito}) turns out to be accurate, not only for large $\nu$ 
but also for small $\nu$. It typically gives
the imaginary part of the zero closest to the imaginary axis with $2-3$ digits for $\nu>4$ and improving as the order 
becomes larger.

\section{{\large Algorithms}}

In \cite{Segura:2012:CCZ}, an algorithm for computing complex zeros of special functions satisfying second order ODEs
$y''(z)+A(z)y(z)=0$ was constructed which was based on the assumption that the coefficient $A(z)$ is ``locally constant".
The idea of the method is to follow the ASLs in the Liouville-Green approximation by taking steps in the direction of the
ASLs, and to apply a fourth order fixed point iteration $x_{n+1}=g(x_n)$ with iteration function
\begin{equation}
\label{FPM4}
g(z)=z-\Frac{1}{\sqrt{A(z)}}\arctan\left(\sqrt{A(z)}\Frac{y(z)}{y'(z)}\right)
\end{equation}
in order to compute the zeros.

Once a first zero $z_1$ is computed by iterating (\ref{FPM4}) with an starting value conveniently
chosen, the method takes a step 
\begin{equation}
\label{steps}
\hat{z}=\pm \pi/\sqrt{A(z_1)}, 
\end{equation}
where the sign is 
preferably chosen in such a way that the step is taken in the direction of decreasing $|A(z)|$. If $A(z)$ was constant,
$\hat{z}$ would be another zero. In general $\hat{z}$ will not be another zero, but it will be a value which gives convergence 
to the next zero in the ASL
by iterating this value with (\ref{FPM4}). Once a second zero $z_2$ is computed by iterating  $\hat{z}$ with
(\ref{FPM4}) we would take another step $\hat{z}=\pm \pi/\sqrt{A(z_2)}$, with the sign chosen as before, and so on.

The main difficulty consists in obtaining good initial values for computing the first
zero in each ASL. Although the good global convergence properties of (\ref{FPM4}) are such that
it is possible to build a method not using such estimations, it is in any case convenient to have sharper first estimates
in order to improve the speed of computation of the first zero. 
On the other hand, it is useful to know in advance which are the approximate
ASLs where the zeros lie in order to avoid testing regions where there are no zeros. This information is contained in the present
paper. Furthermore, if a zero of the solution is known, with this information it is possible to determine the rest of zeros (or
a finite amount of them) reliably.

\subsection{{\normalsize Algorithm for Airy functions}}

Consider first the case of the Airy functions and let suppose that we are interested in computing the zeros of
a solution with a zero at $z_0$. Such solution is proportional to
$y(z)={\rm Bi}(z_0) {\rm Ai}(z)-{\rm Ai}(z_0){\rm Bi}(z_0)$; then, taking 
\begin{equation}
\alpha=\arctan\left(-\Frac{{\rm Ai}(z_0)}{{\rm Bi}(z_0)}\right),
\end{equation}
we have that the zeros of $y(z)$ are the zeros of ${\cal A}(\alpha,z)$.

For computing the zeros of ${\cal A}(\alpha,z)$ we can consider initial values given by the asymptotic approximations, which
will be more accurate as $|z|$ is larger, and then take the steps (\ref{steps}) in the direction of decreasing $|z|$. In fact,
the asymptotic approximations are not really necessary for our method and taking initial values close to the principal 
ASLs is enough
(though using the expansions for estimating a first zero is convenient).

For the zeros of the derivative the same procedure can be considered, but replacing (\ref{FPM4}) by a fixed point giving
convergence to the zeros of  ${\cal A}'(\alpha,z)$. A possibility is to consider
\begin{equation}
\label{FPM2}
\tilde{g}(z)=z+\Frac{1}{\sqrt{A(z)}}\arctan\left(\Frac{1}{\sqrt{A(z)}}\Frac{y'(z)}{y(z)}\right),
\end{equation}
which gives convergence to the zeros of $y'(z)$. Contrary to (\ref{FPM4}), which has order of convergence $4$, 
(\ref{FPM2}) has order $2$, like Newton's method. This fixed point method was discussed in \cite{Gil:2003:CTZ} for the case
of real zeros and it was proved to be globally convergent when $A(z)$ is monotonic \cite[Theorem 3.2]{Segura:2010:RCO}. 
A better possibility consists in applying the idea of \cite{Gil:2012:CTR} for computing zeros of the derivatives of solutions
of ODEs: take the derivative of the differential equation and eliminate $y(z)$ by using the same differential equation; this gives
a second order ODE for $w(z)=y'(z)$; transform to normal form and compute the corresponding fourth order fixed point method.
A straightforward computation gives the following fourth order fixed point method for computing
the zeros of ${\cal A}'(\alpha,z)$:

\begin{equation}
g(z)=z-\Frac{1}{w(z)}\arctan\left(\Frac{w(z)}{z\Frac{{\cal A}(\alpha,z)}{{\cal A}^{\prime}(\alpha,z)}-\Frac{1}{2z}}\right),\,
w(z)=\sqrt{-z-\Frac{3}{4z^2}}.
\end{equation}

\subsection{{\normalsize Algorithm for Bessel functions}}

A simple but effective algorithm for computing complex 
zeros of Bessel functions ${\cal C}_{\nu}(\alpha,z)$ was given in \cite{Segura:2012:CCZ}.
The method can be largely improved by using the information given in the present paper. 
In addition, as in the case of Airy functions,
we can design a method that, given a zero $z_0$ of a solution of the ODE is able to compute the rest of zeros (a finite amount
of them) just by considering that $\alpha=\arctan(J_{\nu}(z_0)/Y_{\nu}(z_0))$. 

As discussed before, there are always zeros as $\Re z\rightarrow +\infty$ except for Hankel functions. We compute these
zeros starting from large $\Re z$ and use (\ref{steps}) with the minus sign. 
We stop when $\Re z<\nu$ is reached (when the region of Airy-type zeros is reached).
The MacMahon expansions can be used for the estimation of a first large zero. A simpler and also very effective
starting value is $z=L-i\Im \alpha$ (see Eqs. (\ref{MMO}) and (\ref{betaeta})), 
with $L$ a positive value that can be chosen at will (depending on how many zeros
are wanted).

With respect to the zeros running parallel to the negative axis as $\Re z\rightarrow -\infty$, the information given in
section \ref{parbra} can be used to determine when there are zeros above and/or below the negative real axis (see also
section \ref{Hankel} for the case of Hankel functions). When there are zeros, it is possible to compute them starting
from large $|z|$ and in the direction of decreasing $|z|$ (using (\ref{steps} with the plus sign) 
until $\Re z>-\nu$; the starting value can be obtained from McMahon-type
expansions, with a first approximation given by (\ref{zsf}) with $m=\pm 1$. It is simple and equally effective to
consider a value $-L-ia$ for zeros above the branch cut and $-L-ib$ below the branch cut, with $a$ and $b$ given in 
(\ref{paraa}) and (\ref{parab}).

Finally, for computing the Airy-type zeros, we can give good starting values using the information of section \ref{ojito}.
We compute the intersection of the ASLs containing these zeros with the imaginary axis solving the equations 
(\ref{masterojito}). After computing one of such intersections, we iterate with (\ref{FPM4}) to compute a first zero;
then we move to the right to compute subsequent zeros using the steps (\ref{steps}) with plus sign; and starting again
from this first zero we move to the left with the minus sign. In both cases, we stop when $|\Re z|>\nu$ is reached or
when the imaginary part of $z$ becomes very small or changes sign (this last criterion is needed for the zeros
inside the eye-shaped region).

For the zeros of the first derivative, the same scheme works, but the fixed point method (\ref{FPM4}) has to be replaced. 
We can use the fixed point iteration developed for the real zeros of
$\gamma {\cal C}_{\nu}(\alpha,z)+z {\cal C}_{\nu}^{\prime}(\alpha,z)$ in \cite{Gil:2012:CTR}, which is of fourth
order (also for complex zeros). Taking $\gamma=0$ the fixed point iteration reads:
\begin{equation}
\label{4DB}
\begin{array}{l}
g(z)=z-\Frac{z(z^2-\nu^2)}{W(z)}\arctan\left(W(z)\Frac{D{\cal C}_{\nu}(\alpha,z)+ E{\cal C}_{\nu+1}(\alpha,z)}
{M{\cal C}_{\nu}(\alpha,z)+ N{\cal C}_{\nu+1}(\alpha,z)} \right),\\
\\
W(z)=\sqrt{z^6+Pz^4 +Q z^2 +R},\\ 
\\
P=-3(\nu^2+1/4),\,Q=\nu^2 (3\nu^2-5/2),\,
R=-\nu^4 (\nu^2-1/4),\\ 
\\
D=\nu,\,E=-z,\,M=-z^4+2\nu (\nu -1/4) z^2-\nu^3 (\nu +1/2),\,N=\frac12 (z^3+\nu^2 z) .
\end{array}
\end{equation}

Maple and Fortran codes implementing these methods are under construction \cite{Gil:2014:ZAM}.

\section*{{\large Acknowledgements}}

The authors acknowledge some financial support from {\emph{Ministerio de Econom\'{\i}a y Competitividad}}, 
project MTM2012-34787. 

\bibliographystyle{plain}
\bibliography{bibliobes}

\end{document}